# Partitioning the set of natural numbers into Mersenne trees and into arithmetic progressions; Natural Matrix and Linnik's constant

Gennady Eremin

**Abstract.** We partition a series of natural numbers into infinite number sequences. We consider two partitioning options: (a) a forest of unary trees with recurrence formula of Mersenne numbers, and (b) a set of arithmetic progressions with difference $2^k$. Every tree starts with an even number, and any even number starts a certain tree. Unary trees with initial terms 0, 2, ..., 22 are in the OEIS. In the partitioning into arithmetic progressions, each progression starts with a Mersenne number, and each Mersenne number is the beginning of a particular arithmetic progression. Unary trees starting from some term are contained in OEIS A036991 (compact Dyck path codes), so we consider A036991 as a backbone of the partitions. In particular, we prove the existence of an arithmetic progression of any length in A036991. As a result of the partitions, we obtain a Natural Matrix with a packing function that captures the bijection between the set of natural numbers and the set of ordered pairs of natural numbers.

In Natural Matrix, the even natural numbers are located on the *x*-axis, so the selection of primes in the considered arithmetic progressions is greatly simplified. A method for proving the infinity of primes in A036991 is proposed. In this regard, an attempt is made to reduce the Linnik's constant to 2.

## 1 Introduction

The natural numbers of the OEIS sequence A036991 (including 0) are Dyck path codes. Dyck paths (aka balanced bracket sequences or Dyck words) are well known in the Catalan family [1]. We called the A036991 terms *Dyck numbers* [2], since in bracket sequences each left parenthesis is encoded by a zero (correspondingly, right parentheses are encoded by a one), and as a result we obtain compact numerical binary codes (in integers leading zeros are usually omitted). In such binary codes, *Dyck dynamics* are monitored: in each binary suffix, the number of zeros must not exceed the number of ones. All Dyck numbers are odd (except 0) because, according to Dyck dynamics, the final binary digit (suffix of length 1) is always 1.

The sequence A036991 includes dozens of sequences from the OEIS [14], among which are the well-known Mersenne numbers (OEIS A000225), which have no zeros in their binary code except for the initial zero. We work with sequences for which the recurrence formula is similar to the Mersenne numbers and has the following form:

(1) $$a_{n+1} = 2a_n + 1 = a_n + (a_n + 1) = a_n + \mathrm{S}a_n \quad (a_0 = 0, 2, 4, \ldots),$$

where S is the successor function for natural numbers (we borrowed the notation from [3]). Further, we denote the series of natural numbers with zero by $\mathbb{N} = \{0, 1, 2, 3, ...\}$ (see OEIS A001477). This notation is the most modern, especially in the field of computer science, where counting often starts from zero.

## 2 Unary Mersenne trees

The recurrent formula (1) works in many OEIS sequences. Formula (1) is linear, and we can say that such a formula realizes an axiomatic approach, something analogous to Peano's Axioms. The axiomatic approach is based on the axiomatization of the properties of ordinal numbers: according to formula (1), every term has a unique *successor* (heir), and every odd term has a unique, possibly even *predecessor* (parent). In this case, the recurrence formula resembles a succession function, and the domain of definition of such a function is a *unary tree* with some even root (it is 0, 2, 4, 6, etc.). The corresponding linear precedence function is $a_{n-1} = (a_n - 1)/2$ for each odd term $a_n$.

In unary trees, it is convenient and logical to work with binary codes of terms; for example, each step of moving along the tree from parent to successor corresponds to adding a unit at the end of the binary code of the parent. Accordingly, when moving backwards from the successor (odd number) to the parent (possibly even number) in the binary code of the successor, we simply remove the final one. In natural numbers, the length of binary codes is not limited, so in each unary tree the number of odd terms is infinite, that is, we can take any number of steps from the root to infinity. The backward movement to the root of the tree is accompanied by a stepwise reduction of the unit suffix in odd terms until the first zero appears at the end of the binary code; such a zero fixes an even number, the root of the given tree. Obviously, we can formulate a corresponding statement.

**Proposition 1.** *Every even natural number, including 0, is a root of a unique unary tree.*

A unary tree with root 0 is the sequence of Mersenne numbers 0, 1, 3, 7, 15, etc. (binary codes are 0, 1, 11, 111, 1111, etc.). When moving backward along such a tree, we can remove the units of the binary code step by step, and as a result we will get a zero root (the leading zero can be prescribed in advance at the beginning of an odd Mersenne number).

The root of the next unary tree is the even number 2, here is the beginning of the corresponding sequence: 2, 5, 11, 23, 47, etc. (binary codes 10, 101, 1011, 10111, 101111, etc.), and this is the OEIS sequence [A153893](A153893). In such a tree, moving back to the root by removing the last unit in the binary codes of the odd numbers at each step, we will eventually return to root 2.

It is not hard to see, such sequences are infinite, because for each term we simply dock a unit step by step at the end of the binary code. The recurrence formula of Mersenne numbers works in unary trees, and it is logical to call such sequences Mersenne trees.

It is not hard to see, such sequences are infinite, because for each term we simply dock a unit at the end of the binary code step by step. The recurrence formula of Mersenne numbers works in unary trees, and it is logical to call such sequences *Mersenne trees*. We denote the Mersenne tree with root 0 by $\mathbb{M}_0$, obviously $\mathbb{M}_0 = $ A000225. The next unary tree with even root 2 is $\mathbb{M}_1 = $ A153893. In general, the Mersenne tree $\mathbb{M}_k$ starts with an even root $2k$. As a result, we obtain an infinite forest of unary trees with even roots, into which the set of natural numbers seems to be split. Later we will consider the corresponding theorem.

Currently, more than 10 unary Mersenne trees are OEIS sequences, indicating their importance. In Table 1, we have shown in the left column the roots of the first unary trees, even numbers from 0 to 32. The second column lists the initial nodes of sequences (odd numbers) that are not A036991 terms; the binary codes of such numbers have suffixes with negative Dyck dynamics. As we can see, the number of such nodes is still small. Degrees of two have the largest number of such terms, for example, for a root of $2^n$ (there are $n$ finite zeros in binary) the first $n$ numbers including the root are not A036991 terms. The third column shows the subsequent terms of the sequence A036991, and there are an infinite number of such terms in each unary tree, since the length of the binary codes is not limited, and any number of ones can be assigned to the end of each code without degrading the Dyck dynamics. We can say that the vast majority of odd nodes of Mersenne trees are chosen from the sequence A036991.

| Roots | Not terms of A036991 | A036991 terms | OEIS |
|---|---|---|---|
| 0 | | 1, 3, 7, 15, 31, 63, 127, 255, 511, … | A000225 |
| 2 | | 5, 11, 23, 47, 95, 191, 383, 767, … | A153893 |
| 4 | 9 | 19, 39, 79, 159, 319, 639, 1279, … | A153894 |
| 6 | | 13, 27, 55, 111, 223, 447, 895, … | A086224 |
| 8 | 17, 35 | 71, 143, 287, 575, 1151, 2303, … | A052996 \{1, 3} |
| 10 | | 21, 43, 87, 175, 351, 703, 1407, … | A086225 |
| 12 | 25 | 51, 103, 207, 415, 831, 1663, … | A198274 |
| 14 | | 29, 59, 119, 239, 479, 959, 1919, … | A196305 |
| 16 | 33, 67, 135 | 271, 543, 1087, 2175, 4351 … | A198275 |
| 18 | 37 | 75, 151, 303, 607, 1215, 2431, … | A198276 |
| 20 | 41 | 83, 167, 335, 671, 1343, 2687, … | A171389 |
| 22 | | 45, 91, 183, 367, 735, 1471, 2943, … | A291557 |
| 24 | 49, 99 | 199, 399, 799, 1599, 3199, 6399, … | |
| 26 | | 53, 107, 215, 431, 863, 1727, … | |
| 28 | 57 | 115, 231, 463, 927, 1855, 3711, … | |
| 30 | | 61, 123, 247, 495, 991, 1983, … | |
| 32 | 65, 131, 263, 527 | 1055, 2111, 4223, 8447, 16895, … | |

Table 1. The first Mersenne trees and corresponding OEIS sequences.

Let's repeat, the set of even natural numbers is infinite, so we get an infinite number of unary Mersenne trees with even roots; such trees do not intersect due to the linearity of the succession formula and, correspondingly, the linearity of the precedence formula. An odd natural number cannot be a node in several unary trees, since we reach a single even number when moving from an odd node to the root. We will give below a rigorous proof of what has been said, but now let's prove the accompanying statement.

Partitioning the set of natural numbers into Mersenne trees and into arithmetic progressions; Natural Matrix and Linnik's constant.

**Proposition 2.** *Every natural number is a node in some Mersenne tree with a particular even root.*

*Proof.* For even natural numbers including 0, we proved the analogous Proposition 1. Let an odd natural number $t$ not be a node in any Mersenne tree, and let the length of the binary code of $t$ be $n$ and the length of the unit suffix be $m$, $n \geq m > 0$. In the case $n = m$ we deal with the Mersenne number $t = M_n$ (there are no zeros in the binary code of $t$), and then the number $t$ is a node of the null Mersenne tree $\mathbb{M}_0$. If $m < n$, it means that in the binary code of $t$ there is a zero at some position $m$ (from the end of $t$), and this zero ends the even number

(2) $\qquad\qquad\qquad r = (t - M_m) / 2^m = \lfloor t / 2^m \rfloor$,

where $\lfloor x \rfloor$ denotes the largest integer less than or equal to $x$, the floor function. As a result, we have an even number $r$ which is the root of the ($r/2$)-th unary Mersenne tree, and there is a node $t$ in this tree at the $m$-th step from the root. We got a contradiction. □

Thus, every natural number including 0 is a node of some Mersenne tree. Let us find a general formula for nodes in such trees. Often the formula for the general term of a sequence is obtained using recurrence relations, for example, here is a well-known formula for Mersenne numbers: $M_n = 2^n - 1$, $n = 0, 1, 2, \ldots$ The last column of Table 1 lists the OEIS sequences where the reader can find similar formulas. Below we give the formulas for the first Mersenne trees.

$\mathbb{M}_0 = \{1 \times 2^n - 1\} = \{0, 1, 3, 7, 15, \ldots\} = $ A000225;
$\mathbb{M}_1 = \{3 \times 2^n - 1\} = \{2, 5, 11, 23, 47, \ldots\} = $ A153893;
$\mathbb{M}_2 = \{5 \times 2^n - 1\} = \{4, 9, 19, 39, 79, \ldots\} = $ A153894;
$\mathbb{M}_3 = \{7 \times 2^n - 1\} = \{6, 13, 27, 55, 111, \ldots\} = $ A086224;
$\mathbb{M}_4 = \{9 \times 2^n - 1\} = \{8, 17, 35, 71, 143, \ldots\} = $ A052996\{1, 3\};
$\mathbb{M}_5 = \{11 \times 2^n - 1\} = \{10, 21, 43, 87, 175, \ldots\} = $ A086225;
………
$\mathbb{M}_{11} = \{23 \times 2^n - 1\} = \{22, 45, 91, 183, 367, \ldots\} = $ A291557.

Obviously, for the $k$th Mersenne tree, $k = 0, 1, 2$, etc., we obtain the following expression:

(3)  $\mathbb{M}_k = \{(2k+1) \times 2^n - 1\} = \{2^{n+1}k + M_n\}$
    $= \{2k, 4k+1, 8k+3, 16k+7, 32k+15, \ldots\}$.

Let us prove the non-intersectionality of unary Mersenne trees using formula (3).

**Proposition 3.** *A natural number cannot be a node in two or more Mersenne trees.*

*Proof.* For even natural numbers, everything is obvious since every unary tree has a unique even root. Let odd numbers $a = (2k+1) \times 2^n - 1$, $n > 0$, and $b = (2l+1) \times 2^m - 1$, $m > 0$, be given, and let $a = b$. In this case $(2k+1) \times 2^n = (2l+1) \times 2^m$, and then $n = m$ (in the decomposition of $a$ and $b$, the number of twos is the same). The result is $2k+1 = 2l+1$, and then $k = l$. Thus we get a single node in a fixed unary tree. □

Now we can formulate the partition theorem for the set of natural numbers (recall that partitioning a set is its representation as a union of an arbitrary number of pairwise disjoint subsets).

**Theorem 4.** *The set of natural numbers including 0 is partitioned into disjoint Mersenne trees with even roots 0, 2, 4, 6, etc. In each such tree the recurrence formula (1) works.*

As a result, we get the following formula:
(4)  $\mathbb{N} = \bigcup_k \mathbb{M}_k$, $k = 0, 1, 2, \ldots$

And at the end of this section, let's talk a little about the generating functions of the described sequences. In combinatorics, generating functions allow us to work with different objects using analytical methods. For example, analytic functions are often used to obtain an explicit formula for the general term of a sequence. Typically, the generating function is a formal power series. For example, the generating function of the Mersenne numbers (0th Mersenne tree) is

(5) $\quad G_0(x) = 0x^0 + 1x^1 + 3x^2 + 7x^3 + 15x^4 + \ldots + M_n x^n + \ldots = x/(1 - 3x + 2x^2).$

The last equality in (5) is easy to obtain (see Example 5 in [4]). Let's borrow generating functions for other Mersenne trees from some OEIS sequences.

1st Mersenne tree, A153893: $G_1(x) = (2 - x) / (1 - 3x + 2x^2).$
2nd Mersenne tree, A153894: $G_2(x) = (4 - 3x) / (1 - 3x + 2x^2).$
3rd Mersenne tree, A086224: $G_3(x) = (6 - 5x) / (1 - 3x + 2x^2).$
………
11th Mersenne tree, A291557: $G_{11}(x) = (22 - 21x) / (1 - 3x + 2x^2).$

Obviously, the generating function for the kth Mersenne tree is

(6) $\quad G_k(x) = (2k - (2k-1)x) / (1 - 3x + 2x^2).$

## 3  Natural Matrix

Let's return to Table 1 and work a little with formula (3). In the first column of Table 1 we have initial terms of unary trees 0, 2, 4, ... (even natural numbers with zero, the OEIS sequence A005843); i.e. we have an arithmetic progression with initial term $a_0 = 0$ and difference 2 ($a_1 = a_0 + 2 = 2$, etc.). Let us choose an initial term from formula (3) (option $n = 0$), and write the first progression as $A_0(k) = \{M_0 + 2^{0+1}k\} = \{2k\}$, recall in this case $k$ is the index of the Mersenne tree.

The following terms in unary trees (option $n = 1$) are another arithmetic progression 1, 5, 9, ... (the OEIS sequence A016813) with initial term $a_0 = M_1 = 1$ and a difference 4 ($a_1 = a_0 + 4 = 5$, etc.). Let us write the following arithmetic progression in the form $A_1(k) = \{M_1 + 2^{1+1}k\} = \{1 + 4k\} = $ A016813. This procedure can be continued infinitely. Let's write down the first arithmetic progressions and the corresponding generating functions.

$A_0(k) = \{0 + 2k\} = \{0, 2, 4, 6, 8, 10, \ldots\} = $ A005843;   $G_0(x) = 2x/(1-x)^2.$
$A_1(k) = \{1 + 4k\} = \{1, 5, 9, 13, 17, 21, \ldots\} = $ A016813;   $G_1(x) = (1+3x)/(1-x)^2.$
$A_2(k) = \{3 + 8k\} = \{3, 11, 19, 27, 35, 43, \ldots\} = $ A017101;   $G_2(x) = (3+5x)/(1-x)^2.$
$A_3(k) = \{7 + 16k\} = \{7, 23, 39, 55, 71, 87, 103, 119, \ldots\},$   $G_3(x) = (7+9x)/(1-x)^2.$
$A_4(k) = \{15 + 32k\} = \{15, 47, 79, 111, 143, 175, 207, \ldots\},$   $G_4(x) = (15+17x)/(1-x)^2.$

We borrowed the generating functions $G_0(x)$, $G_1(x)$ and $G_2(x)$ from the corresponding OEIS sequences, other generating functions are not difficult to obtain, the reader can easily check (see [4]). The general form of the $n$th arithmetic progression and its derivative function are as follows:

(7) $\quad A_n(k) = \{M_n + 2^{n+1}k\} = \{M_n, M_n + 1 \times 2^{n+1}, M_n + 2 \times 2^{n+1}, M_n + 3 \times 2^{n+1}, \ldots\},$

$$G_n(x) \;=\; (M_n + (M_n+2)\,x)/(1-x)^2, \text{ for } n = 0, 1, 2, \ldots$$

Thus we can talk about some infinite rectangular matrix of natural numbers, let's call it a *Natural Matrix* and denote it by $\mathcal{M}$. The rows of the matrix are Mersenne trees, and the columns are arithmetic progressions, with each successive column doubling the difference in progressions. In the following table we have shown Narural Matrix with the initial element $m_{0,0} = 0$ (upper left corner of the matrix); in the lines with dark fill are marked initial numbers, which are not A036991 terms.

|    | 0  | 1  | 2   | 3   | 4   | 5    | 6    | 7    | 8    | 9     | 10    | … |
|----|----|----|-----|-----|-----|------|------|------|------|-------|-------|---|
| 0  | 0  | 1  | 3   | 7   | 15  | 31   | 63   | 127  | 255  | 511   | 1023  | … |
| 1  | 2  | 5  | 11  | 23  | 47  | 95   | 191  | 383  | 767  | 1535  | 3071  | … |
| 2  | 4  | 9  | 19  | 39  | 79  | 159  | 319  | 639  | 1279 | 2559  | 5119  | … |
| 3  | 6  | 13 | 27  | 55  | 111 | 223  | 447  | 895  | 1791 | 3583  | 7167  | … |
| 4  | 8  | 17 | 35  | 71  | 143 | 287  | 575  | 1151 | 2303 | 4607  | 9215  | … |
| 5  | 10 | 21 | 43  | 87  | 175 | 351  | 703  | 1407 | 2815 | 5631  | 11263 | … |
| 6  | 12 | 25 | 51  | 103 | 207 | 415  | 831  | 1663 | 3327 | 6655  | 13311 | … |
| 7  | 14 | 29 | 59  | 119 | 239 | 479  | 959  | 1919 | 3839 | 7679  | 15359 | … |
| 8  | 16 | 33 | 67  | 135 | 271 | 543  | 1087 | 2175 | 4351 | 8703  | 17407 | … |
| 9  | 18 | 37 | 75  | 151 | 303 | 607  | 1215 | 2431 | 4863 | 9727  | 19455 | … |
| 10 | 20 | 41 | 83  | 167 | 335 | 671  | 1343 | 2687 | 5375 | 10751 | 21503 | … |
| 11 | 22 | 45 | 91  | 183 | 367 | 735  | 1471 | 2943 | 5887 | 11775 | 23551 | … |
| 12 | 24 | 49 | 99  | 199 | 399 | 799  | 1599 | 3199 | 6399 | 12799 | 25599 | … |
| 13 | 26 | 53 | 107 | 215 | 431 | 863  | 1727 | 3455 | 6911 | 13823 | 27627 | … |
| 14 | 28 | 57 | 115 | 231 | 463 | 927  | 1855 | 3711 | 7423 | 14847 | 29695 | … |
| 15 | 30 | 61 | 123 | 247 | 495 | 991  | 1983 | 3967 | 7935 | 15871 | 31743 | … |
| 16 | 32 | 65 | 131 | 263 | 527 | 1055 | 2111 | 4223 | 8447 | 16895 | 33791 | … |
| 17 | 34 | 69 | 139 | 279 | 559 | 1119 | 2239 | 4479 | 8959 | 17919 | 35839 | … |
| 18 | 36 | 73 | 147 | 295 | 591 | 1183 | 2367 | 4735 | 9471 | 18943 | 37887 | … |
| …  | …  | …  | …   | …   | …   | …    | …    | …    | …    | …     | …     | … |

Table 2. Infinite Natural Matrix $\mathcal{M}$.

As we known, the asymptotic density of the set of even natural numbers $A_0(k)$ is equal to ½, i.e. even numbers make up half of the natural series (the same is the density of odd natural numbers); the density of the next arithmetic progression $A_1(k)$ is ¼ (the quarter of natural numbers, the difference of the progression is 4) and so on. Like Mersenne trees the arithmetic progressions also have no common terms, so the union of such progressions gives us the set of natural numbers, since

$$\tfrac{1}{2} + \tfrac{1}{4} + \tfrac{1}{8} + \ldots \;=\; 1.$$

Let's formulate the corresponding theorem.

**Theorem 5.** *The union of arithmetic progressions $A_n(k)$, $n = 0, 1, 2$, etc, is a series of natural numbers.*

(8) $$\mathbb{N} \;=\; \bigcup_n A_n(k), \; n = 0, 1, 2, \ldots$$

As we see, in Natural Matrix $\mathcal{M}$, both row union and column union gives us the set $\mathbb{N}$. Naturally, the union of columns and the union of rows of the same matrix must coincide, so the proof of Theorem 5 confirms Theorem 4 and vice versa. For the terms of the matrix $\mathcal{M} = \{m_{i,j}\}$, $i, j = 0, 1, 2$, etc., it is not difficult to obtain some formulas, for example, $m_{0,j} = M_j$, $m_{i,0} = 2i$, $m_{i,i} = M_i + 2^{i+1}i$.

Let us formulate the last statement in this section, the statement of bijection between the sets $\mathbb{N}$ and $\mathcal{M}$.

**Proposition 6.** *The mapping between the natural numbers, including 0, and the elements of Natural Matrix $\mathcal{M}$ is a bijection.*

*Proof.* According to formula (3), the matrix element $m_{i,j}$ corresponds to a natural number, namely: $m_{i,j} = 2^{j+1}i + M_j = (2i+1) \times 2^j - 1 \in \mathbb{N}$. In the case of formula (3), $k = i$ is a line in $\mathcal{M}$ and $n = j$ is a column in $\mathcal{M}$. Conversely, for an arbitrary natural number $n \in \mathbb{N}$ let's find the corresponding element $m_{i,j} \in \mathcal{M}$. There are three options:

  (i) $n$ is an even number, then $n = m_{i,0}$, where $i = n/2$;
  (ii) $n$ is a $k$-th Mersenne number, then $n = m_{0,k}$;
  (iii) in the binary code of number $n$, the length of the unit suffix is $k$, and then $n = m_{i,k}$, where $i = (n - M_k)/2^{k+1}$ (see the proof of Proposition 2).

The assertion is proven. □

## 4 Packing function of Natural Matrix

A bijection between the set of natural numbers $\mathbb{N}$ and a natural matrix $\mathcal{M}$ can be viewed as a self-inverse permutation of non-negative integers. The one-to-one correspondence between the set $\mathbb{N}$ and the set of ordered pairs $\mathbb{N} \times \mathbb{N}$ is known in the literature, and it is the *packing* (or *pairing*) *function*, the bijective mapping $\pi : \mathbb{N} \times \mathbb{N} \to \mathbb{N}$. Symmetric *Cantor polynomials* are well known [5, 6, 7]:

(9) $\quad \pi(x, y) = (3x + (x+y)^2 + y)/2 \quad$ and $\quad \pi'(x, y) = (x + (x+y)^2 + 3y)/2 = \pi(y, x)$.

Currently, mathematicians believe that no other polynomial maps bijectively $\mathbb{N} \times \mathbb{N}$ onto $\mathbb{N}$. This is exactly what the famous theorem of R. Fueter and G. Pólya [8] states: "*There are only two quadratic packing polynomials, and they are Cantor polynomials*".

In our case (a) from the formula for the general term of the $k$-th Mersenne tree $(2k+1) \times 2^n - 1$, see (3), and (b) from the formula for the general term of the $n$-th arithmetic progression $2^{n+1}k + M_n$ with difference $2^{n+1}$, see (7), we obtain the following symmetric packing functions:

(10) $\quad\quad F(x, y) = (2x+1) \times 2^y - 1 \quad$ and $\quad G(x, y) = (2y+1) \times 2^x - 1 = F(y, x);$

and these are not Cantor polynomials. Functions (10), like Cantor polynomials, take integer values at points of the plane with non-negative integer coordinates and map each pair of such non-negative numbers into one unique non-negative integer. We can say that we get lossless compression of numerical data.

The packing functions (10) are well known; it is verified that they are bijective and monotone for each argument (see [9, page 47] and [10, page 7, Figure 3]). Functions are not difficult to reverse if we decompose, for example, $F(x, y) +1$ by the degree of two, $2^y$, and an odd factor $2x +1$. Such functions are important in logic and recursive function theory. Functions (10) group even natural numbers (in the case of the matrix $\mathcal{M}$ this is the initial column, see Table 2), and then, for example, the search for prims is greatly simplified, because we know and easily accessible for analysis a compact part of the data array, which contains only odd natural numbers.

Next, we will work with the first packing function $F(x, y)$. Below in Figure 1 we have shown the initial rows of the function $F(x, y)$ for the matrix $\mathcal{M}$. The numbers that are not terms of A036991 are shown in small font and pale. In the rows, we connected neighboring A036991 terms with dots.

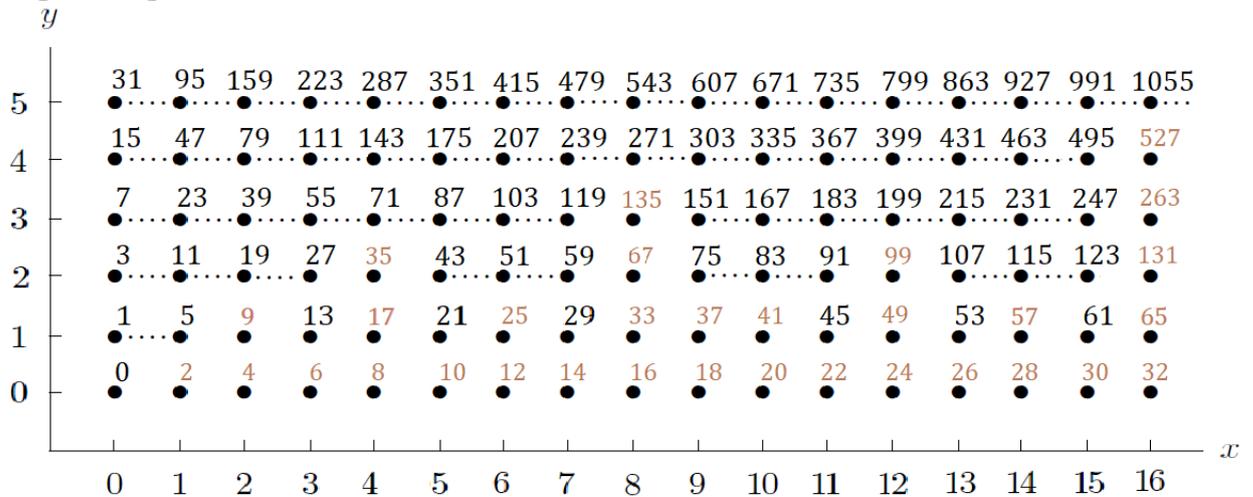

Figure 1. The packing function $F(x, y)$.

In Figure 1, in line $y$, we have an arithmetic progression with the initial term $M_y$ and the difference $2^{y+1}$. Further in the lines we will be interested only in the initial chains of Dyck numbers, initial segments (intervals) with A036991 terms. In the lowest row ($y = 0$) we have even natural numbers, the difference of such a progression $d_0 = 2$, and only the first integer 0 is the A036991 term; let's denote such a segment $S_0 = \{0\}$, the length of the segment is $\#S_0 = 1$. In the next line ($y = 1$) we have the arithmetic progression A016813 with two Dyck numbers in the first segment, they are 1 and 5, that is, $S_1 = \{1, 5\}$, $\#S_1 = 2$; the difference of the progression is also doubled, $d_1 = 4$. Then on the line $y = 2$ we get a segment with four Dyck numbers $S_2 = \{3, 11, 19, 27\}$, $\#S_2 = 4$ and $d_2 = 8$. The trend is obvious, at the beginning of each subsequent line the number of A036991 terms is doubled, and the difference of the corresponding arithmetic progression is also doubled.

Let's demonstrate on binary codes the generation of initial segments in matrix rows. In the Mersenne number $M_y$ (remember, in the binary expansion of such a number there are no zeros, but only $y$ units), we add the binary prefix 10 at the beginning of the code; this is equivalent to adding the integer $d_y = 2^{y+1}$ (the difference of the corresponding arithmetic progression). And then let's add $d_y$ step by step until we get the same Mersenne number in the prefix. Below we have shown the generation of segment $S_3$ for example. The red color shows the incremental prefix: 0, 1, 10, ..., 111.

<span style="color:red">00</span>111, <span style="color:red">10</span>111, <span style="color:red">100</span>111, <span style="color:red">110</span>111, <span style="color:red">1000</span>111, <span style="color:red">1010</span>111, <span style="color:red">1100</span>111, <span style="color:red">1110</span>111.

The first term in $S_3$ is the Mersenne number $M_3 = 111_2 = 7$, the last term is the number $M_3 \times d_3 + M_3 = 7 \times 16 + 7 = 1110111_2 = 119$. The length of $S_3$ is $M_3 + 1 = d_3/2 = 16/2 = 8$.

Thus, at the beginning of line $y$ we get the following segment:

the first term is the $y$-th Mersenne number $M_y = 2^y - 1$,
the last term is $M_y(d_y + 1)$,
the number of Dyck numbers in such a segment is $\#S_y = M_y + 1 = 2^y$,
the difference of the arithmetic progression in the line $y$ is $d_y = 2 \times \#S_y = 2^{y+1}$.

Obviously, it is not difficult to reverse the packing function $F(x, y)$. For example, for an arbitrary natural number $n \in \mathbb{N}$ with a unit binary suffix of length $k$, we obtain in Figure 1 the matrix term $m_{x,y} \in \mathcal{M}$ with the following coordinates:

(11) $\qquad\qquad y = k \quad \text{and} \quad 2x = (n+1)/2^m - 1.$

Let's write down a few initial segments from Figure 1.

$S_0 = \{0\}$, $\#S_0 = 1$, $d_0 = 2$;
$S_1 = \{1, 5\}$, $\#S_1 = 2$, $d_1 = 4$;
$S_2 = \{3, 11, 19, 27\}$, $\#S_2 = 4$, $d_2 = 8$;
$S_3 = \{7, 23, 39, 55, 71, 87, 103, 119\}$, $\#S_3 = 8$, $d_3 = 16$;
$S_4 = \{15, 47, 79, 111, 143, 175, 207, 239, 271, 303, \ldots, 495\}$, $\#S_4 = 16$, $d_4 = 32$;
$S_5 = \{31, 95, 159, 223, 287, 351, 415, 479, 543, 607, \ldots, 2015\}$, $\#S_5 = 32$, $d_5 = 64$;
$S_6 = \{63, 191, 319, 447, 575, 703, 831, 959, 1087, 1215, \ldots, 8127\}$, $\#S_6 = 64$, $d_6 = 128$;
$S_7 = \{127, 383, 639, 895, 1151, 1407, 1663, 1919, 2175, \ldots, 32639\}$, $\#S_7 = 128$, $d_7 = 256$.

Thus, at the beginning of the $y$-th line, $y > 0$, we get a shortened arithmetic progression $S_y$ with difference $2^{y+1}$, which contains $2^y$ Dyck numbers. This is the progression:

(12) $\quad S_y = \{\, M_y = 2^y - 1,\ 3 \times 2^y - 1,\ 5 \times 2^y - 1,\ 7 \times 2^y - 1,\ \ldots,\ M_{y+1} \times 2^y - 1 = M_y(d_y + 1)\,\}.$

Using equality (12), it is easy to find an arithmetic progression of any length in A036991. For example, for arbitrary $k > 1$, the segment $S_l$, $l = \lceil \log_2 k \rceil$, has length $l \geq k$ ($\lceil x \rceil$ denotes the least integer more than or equal to $x$, the ceiling function); and in such a segment we can always choose an arithmetic progression of length $k$. Let's formulate a corresponding statement.

**Theorem 7.** *The OEIS sequence A036991 includes an arithmetic progression of any length.*

In segment (12), the last number is a multiple of $M_y$ and therefore composite for $y > 1$. These numbers are well known, the OEIS has the sequence [A129868](A129868) (the binary codes are [A138148](A138148)), which contains all the maximum terms of (12). These terms are:

$\qquad$ 0, 5, 27, 119, 495, 2015, 8127, 32639, 130815, 523775, 2096127, 8386559, …

The sequence A129868 is contained entirely in A036991.

# 5 Prime numbers in $S_y$, Linnik's constant.

In this paper, we obtained the following two partitions of the set of natural numbers: (a) the forest of unary Mersenne trees and (b) the set of disjoint arithmetic progressions. Both partitions are closely related to the OEIS sequence A036991. The resulting constructions may help to deal with problems that the author has encountered recently, and one of them is the proof of the infinity of primes in A036991. There are no even natural numbers in the sequence A036991, so there are a lot of primes, and this is visible at the beginning of the sequence. In the OEIS there is also a sequence A350577, which collects all prime terms from A036991 and in which the author hypothesized the infinity of the set of primes in A036991/A350577.

In mathematics, the *Dirichlet's theorem* is known about the infinity of primes in an arithmetic progression, provided that the first term of the progression $a$ and the difference $d$ are coprime integers. In the arithmetic progressions we have considered (see Figure 1), this condition is satisfied because every $y$-th progression, $y > 0$, starts with an odd number $a_y = 2^y - 1$, and the difference of the progression is always the power of two, $d_y = 2^{y+1}$ (in our progressions, the difference is twice the first number). There are an infinite number of such arithmetic progressions, and to prove the infinity of primes in A036991, we need only show that each progression has at least one Dyck prime number.

In our case, each arithmetic progression contains both terms and non-terms of A036991, so it is more convenient for us to work with segments (12), which contain only A036991 terms. However, the number of numbers in such segments is limited, despite the fact that in each subsequent segment the length is doubled, i.e., $\#S_{y+1} = 2 \times \#S_y$. We need certain guarantees that starting from some very large segment, we will get at least one prime Dyck number in each subsequent segment. And here let's use well-known Linnik's theorem, which is an amplification of Dirichlet's theorem [11, 12]. Linnik's theorem gives us an upper bound on the value of primes whose existence is proved by Dirichlet's theorem. First, let's give the definition that is used in Linnik's theorem.

**Definition 8.** *Denote by $p(a, d)$ the least prime number in an arithmetic progression of the form $a + nd$, $n \in \mathbb{N}$.*

Let's now formulate Linnik's theorem on prime numbers in arithmetic progression.

**Theorem 9 (Linnik, 1944).** *There exist positive constants $C$ and $L$ such that for any coprime $a, d$ ($1 \leq a < d$) the value of $p(a, d)$ does not exceed $Cd^L$.*

In Theorem 9, the constant $L$ is known as the *Linnik's Constant*. Over the years the Linnik's constant has been redused many times, the last value $L = 5$ was obtained by the German mathematician T. Xylouris in 2011 [13]. Additionally, we note that some authors currently adhere to the hypothesis that $L = 2$. Note that Linnik's theorem introduces an additional constraint $a < d$ (and it is not in Dirichlet's theorem); we have this constraint satisfied in every segment $S_y$, recall that the difference of the arithmetic progression is twice the first number, $d_y = 2(M_y + 1)$.

Let's try to estimate the Linnik's constant for our arithmetic progressions. In particular, we will determine whether the prime number of the segment (if any) corresponds to the current value $L = 5$. In segment $S_y$, the maximum term is

(13) $\quad \max S_y = M_{y+1} \times 2^y - 1 = M_y(2M_y+3) = d_y(d_y-1)/2 - 1 < \tfrac{1}{2} d_y^2,$

and this is significantly less than the current Linnik's constant (reminiscent of the hypothesis $L = 2$). Obviously, if we prove that in the $y$-th arithmetic progression, $y > 0$, there is a prime number $p_{min}(y) < \max S_y$, then we get a decrease of the Linnik's constant, namely,

(14) $\quad p_{min}(y) = p(M_y, d_y) < \tfrac{1}{2} d_y^2 \quad$ for any $y > 0.$

Below in Table 3, we have shown the number of primes in $S_y$, starting from $y = 1$. The primes in the segments are marked in red color. As we can see, there are many primes in the segments, and in each subsequent segment the number of primes increases significantly (almost 2 times). The percentage of primes is quite high, although it gradually decreases as $y$ increases.

| $y$ | $S_y$ | $d_y$ | #primes / #$S_y$ | % of primes |
|---|---|---|---|---|
| 1 | 1, 5 | 4 | 1 / 2 | 50.000 |
| 2 | 3, 11, 19, 27 | 8 | 3 / 4 | 75.000 |
| 3 | 7, 23, 39, 55, 71, 87, 103, 119 | 16 | 4 / 8 | 50.000 |
| 4 | 15, 47, 79, 111, 143, …, 495 | 32 | 7 / 16 | 43.750 |
| 5 | 31, 95, 159, 223, …, 2015 | 64 | 11 / 32 | 34.375 |
| 6 | 63, 191, 319, 447, …, 8127 | 128 | 13 / 64 | 20.312 |
| 7 | 127, 383, 639, … , 32639 | 256 | 24 / 128 | 18.750 |
| 8 | 255, 767, 1279 …, 130815 | 512 | 52 / 256 | 20.312 |
| 9 | 511, 1535, 2559, …, 523775 | 1024 | 95 / 512 | 18.554 |
| 10 | 1023, 3071, …, 2096127 | 2048 | 145 / 1024 | 14.174 |
| 11 | 2047, 6343, …, 8386559 | 4096 | 275 / 2048 | 13.434 |
| 12 | 4095, 12287, …, 33550335 | 8192 | 503 / 4096 | 12.283 |
| 13 | 8191, 24575, …, 134209535 | 16384 | 921 / 8192 | 11.244 |
| 14 | 16383, 49151, …, 536854527 | 32768 | 1717 / 16384 | 10.480 |
| 15 | 32767, 98303, …, 2147450879 | 65536 | 3151 / 32768 | 9.616 |
| 16 | 65535, 196607, …, 8589869055 | 131072 | 5960 / 65536 | 9.094 |
| 17 | 131071, …, 34359607295 | 262144 | 11188 / 131072 | 8.536 |
| 18 | 262143, …, 137438691327 | 524288 | 21171 / 262144 | 8.076 |
| 19 | 524287, …, 549755289599 | 1048576 | 40342 / 524288 | 7.694 |
| 20 | 1048575, …, 2199022206975 | 2097152 | 76511 / 1048576 | 7.296 |
| 21 | 2097151, …, 8796090925055 | 4194304 | 145706 / 2097152 | 6.948 |
| 22 | 4194303, …, 35184367894527 | 8388608 | 277822 / 4194304 | 6.624 |

Table 3. Prime numbers in the initial segments.

Recall that to prove the infinity of primes in A036991, it is enough for us to show the presence of one prime number in each segment $S_y$, starting with some even very large $y$. And by the way, in Table 3, the last 22nd short segment contains more than 4 million Dyck numbers, including more than a quarter of a million primes. To estimate the value of the Linnik's constant for our arithmetic progressions, let's clarify where the least prime number appears in the segments. In the following table, we give in the segment $S_y$, $y > 0$, the position of the first prime number, i.e., its $x$-coordinate (0 corresponds to the initial number in the arithmetic progression, and this is a Mersenne prime).

| y | x | y | x | y | x | y | x | y | x | y | x | y | x | y | x | y | x |
|---|---|---|---|---|---|---|---|---|---|---|---|---|---|---|---|---|---|
| 1 | 1 | 16 | 8 | 31 | 0 | 46 | 11 | 61 | 0 | 76 | 1 | 91 | 16 | 106 | 16 | 121 | 13 |
| 2 | 0 | 17 | 0 | 32 | 2 | 47 | 15 | 62 | 26 | 77 | 27 | 92 | 37 | 107 | 0 | 122 | 13 |
| 3 | 0 | 18 | 1 | 33 | 15 | 48 | 2 | 63 | 4 | 78 | 32 | 93 | 22 | 108 | 61 | 123 | 61 |
| 4 | 1 | 19 | 0 | 34 | 1 | 49 | 9 | 64 | 1 | 79 | 34 | 94 | 1 | 109 | 4 | 124 | 8 |
| 5 | 0 | 20 | 8 | 35 | 12 | 50 | 5 | 65 | 33 | 80 | 7 | 95 | 24 | 110 | 125 | 125 | 7 |
| 6 | 1 | 21 | 3 | 36 | 8 | 51 | 10 | 66 | 31 | 81 | 24 | 96 | 8 | 111 | 37 | 126 | 5 |
| 7 | 0 | 22 | 16 | 37 | 10 | 52 | 32 | 67 | 21 | 82 | 7 | 97 | 25 | 112 | 98 | 127 | 0 |
| 8 | 2 | 23 | 6 | 38 | 1 | 53 | 73 | 68 | 31 | 83 | 45 | 98 | 31 | 113 | 42 | 128 | 47 |
| 9 | 3 | 24 | 19 | 39 | 12 | 54 | 2 | 69 | 25 | 84 | 44 | 99 | 4 | 114 | 17 | 129 | 57 |
| 10 | 2 | 25 | 28 | 40 | 53 | 55 | 1 | 70 | 13 | 85 | 22 | 100 | 38 | 115 | 15 | 130 | 133 |
| 11 | 1 | 26 | 5 | 41 | 7 | 56 | 16 | 71 | 36 | 86 | 127 | 101 | 40 | 116 | 7 | 131 | 27 |
| 12 | 2 | 27 | 10 | 42 | 16 | 57 | 25 | 72 | 2 | 87 | 96 | 101 | 16 | 117 | 22 | 132 | 53 |
| 13 | 0 | 28 | 13 | 43 | 1 | 58 | 38 | 73 | 7 | 88 | 64 | 103 | 1 | 118 | 115 | 133 | 9 |
| 14 | 2 | 29 | 3 | 44 | 17 | 59 | 22 | 74 | 10 | 89 | 0 | 104 | 86 | 119 | 22 | 134 | 5 |
| 15 | 4 | 30 | 106 | 45 | 3 | 60 | 8 | 75 | 12 | 90 | 53 | 105 | 12 | 120 | 76 | 135 | 96 |

Table 4. Position of the first prime number in $S_y$, $x$-coordinate (0 corresponds to $M_y$).

All 135 processed segments have primes, and it is very important that in each segment the prime number appears almost at the very beginning of the segment. As we see, the least prime number appears in the first, maximum in the second hundred terms, and this occurs in those "short" segments whose length is measured in tens of decimal digits. In Table 4, there are only 5 cases where the $x$-coordinate of the least prime number exceeds 100, and these occur in segments 30, 86, 110, 118, and 130. Let's show what we have in segment $S_{130}$:

initial term $\quad a = M_{130} \quad = 1361129467683753853853498429727072845823$,
segment length $\#S_{130} = a + 1 = 1361129467683753853853498429727072845824$,
difference $\quad d = 2 \times \#S_{130} = 2722258935367507707706996859454145691648$,
least prime term $\quad a + 133d = 363421567871562278978884080737128449835007$.

Obviously, in the following segments, it is logical to expect larger $x$-coordinates for the first primes, but the trend is obvious: in each subsequent segment the length is doubled, so the relative position of the least prime number in the segment compared to the length of the entire segment rapidly decreases and approaches the $x$-coordinate of the first term of the progression, i.e. to 0. In other words, the ratio of the position of the first prime number in the segment $S_y$ to the segment length approaches zero with increasing $y$. In the 130th segment discussed above, this ratio is

$$133 / \#S_{130} = 133 / 1361129467683753853853498429727072845824 < 10^{-37}.$$

Such almost zero positions of the least primes allow us to speculate about lower values of the Linnik's constant for our arithmetic progressions. Accordingly, there appear prospects for a constructive proof of the infinity of primes in the sequence A036991.

Finally, let us formulate a hypothesis for the Linnik constant with respect to our arithmetic progressions, this hypothesis is based on two formulas (13) and (14).

**Conjecture 9.** *For arithmetic progressions (7), the Linnik's constant is less than 2, or*

$$p(a, d) < \tfrac{1}{2} d^2.$$

The coefficient $C = \tfrac{1}{2}$ can be considered a consequence of clearing our arithmetic progressions from even numbers. In the case of the proof of Conjecture 9 (a constructive proof is planned in the next paper), the following theorem will become obvious.

**Theorem 10.** *There are an infinite number of primes in the OEIS sequence A036991.*

## Acknowledgements.
The author would like to thank Olena G. Kachko (Kharkiv National University of Radio Electronics) for discussion of the considered integer sequences. Also the author would like to thank the anonymous referees for their valuable comments and suggestions.

## References

[1] R. Stanley. *Catalan numbers*. Cambridge University Press, Cambridge, 2015.

[2] G. Eremin. *Dyck numbers, I. Successor function*, 2022.
arXiv preprint arXiv:2210.00744

[3] NaturalNumbers, www.cs.yale.edu/homes/aspnes/pinewiki/NaturalNumbers.html

[4] Math 3214, Notes on generating functions, 2012.
http://www.math.ualberta.ca/~isaac/math324/s12/gen_functions.pdf

[5] G. Cantor. *Ein beitrag zur mannigfaltigkeitslehre*. Journal für die reine und angewandte Mathematik, 84:242–258, 1878.

[6] M. B. Nathanson. *Cantor polynomials and the Fueter-Pólya theorem*.
arXiv preprint arXiv:1512.08261

[7] B. Sury and M. Vsemirnov. *Packing polynomials on irrational sectors*. Research in Number Theory, 8(3):39, 2022. https://www.isibang.ac.in/~sury/packing.pdf

[8] R. Fueter and G. Pólya. *Rationale abzählung der gitterpunkte*. Vierteljschr. Naturforsch. Ges. Zürich, 58:380–386, 1923.

[9] M. D. Davis, E. J. Weyuker. *Computability, complexity, and languages*. Fundamentals of Theoretical Computer Science. ACADEMIC PRESS, New York 1983.

[10] Matthew P. Szudzik. *The Rosenberg-Strong pairing function*, 2019.
arXiv preprint arXiv:1706.04129

[11] Yu. V. Linnik. *On the least prime in an arithmetic progression I. The basic theorem*. In: Rec. Math. [Mat. Sbornik] N.S. 15(57), 139–178, 1944.

[12] Eric W. Weisstein. *Linnik's Theorem*. From MathWorld – A Wolfram Web Resource. https://mathworld.wolfram.com/LinniksTheorem.html

[13] T. Xylouris. *On the least prime in an arithmetic progression and estimates for the zeros of Dirichlet L-functions*. Acta Arith., 2011, 150: 65–91.

[14] OEIS Foundation Inc., *The On-Line Encyclopedia of Integer Sequences*. Published electronically at http://oeis.org, 2024

Concerned with sequences: A000225, A001477, A005843, A016813, A017101, A036991, A052996, A061579, A086224, A086225, A129868, A138148, A153893, A153894, A171389, A196305, A198274, A198275, A198276, A291557, and A350577.

*Email address*: ergenns@gmail.com
Written: May 26, 2024